\documentclass[12pt]{article}
\usepackage{amssymb}

\usepackage{amsmath}

\usepackage[dvips]{epsfig,graphics}

\newcommand{\A}{\mbox{${\cal A}$}}
\newcommand{\B}{\mbox{$B$}}
\newcommand{\R}{\mbox{$R$}}
\newcommand{\cR}{\mbox{$R$}}

\newcommand\Z{{\mathbb{Z}}}
\newcommand{\Ai}{\mbox{$\mathbb{Z}[x]/(x^2)$}}
\newcommand{\Atwo}{\mbox{$\mathbb{Z}[x]/(x^2)$}}
\newcommand{\Athree}{\mbox{$\mathbb{Z}[x]/(x^3)$}}
\newtheorem{theorem}{Theorem}[section]

\newtheorem{corollary}[theorem]{Corollary}

\newtheorem{definition}[theorem]{Definition}
\newtheorem{example}[theorem]{Example}

\newtheorem{lemma}[theorem]{Lemma}

\newtheorem{proposition}[theorem]{Proposition}
\newtheorem{remark}[theorem]{Remark}

\newtheorem{assumptions}[theorem]{Assumptions}
\newenvironment{proof}[1][Proof]{\textbf{#1.} }{\ \rule{0.5em}{0.5em}}
\input {epsf}
\input{psfig.sty}

\begin{document}

\title{Graph Cohomologies from Arbitrary Algebras}
\author{Laure Helme-Guizon and Yongwu Rong}
\maketitle

\begin{abstract}
\noindent For each commutative, graded algebra with  finite
dimension in each degree, we construct a graded cohomology theory
for graphs whose graded Euler characteristic is the chromatic
polynomial of the graph. This extends our previous work which was
based on the algebra $\mathbb {Z}[x]/(x^2)$.
\end{abstract}

 \tableofcontents

\section{Introduction}

In \cite{K}, Khovanov introduced a graded cohomology theory for
classical links and showed it yields the Jones polynomial by taking
the graded Euler characteristic. This construction has sparked a
good deal of interests in recent years. In \cite{HR}, a graded
cohomology theory for graphs was constructed. The graded Euler
characteristic of the cohomology groups is the chromatic polynomial
of the graph.

Both constructions in \cite{K} and in \cite{HR} depend on a given
graded algebra (more precisely a Frobenius algebra in the case of
knots and links) which is the building block for the chain groups
in the chain complex. However, the amounts of choices of algebras
are quite different. In the case of links, the choices are quite
limited due to invariance of Redemeister moves. For the case of
graphs, the choices are abundant. The algebra used in \cite{HR},
$\Ai$, is the simplest natural choice. The purpose of this note is
to show the construction in \cite{HR} can be made for any graded
$R$-algebra that is finite dimensional in each degree,
 commutative, and
whose product is degree preserving, i.e. $\deg (yz)=\deg y + \deg
z$ for all homogeneous elements $y$ and $z$, where $\R$ is an
integral domain.

In section 2, we explain the definition of the chain complex, and
show that the Euler characteristic of the cohomology groups is
equal to the chromatic polynomial of the graph evaluated at
$\lambda = q\dim \A $ where $q \dim \A $ is the graded dimension
of the algebra $\A$. In section 3, we discuss some basic
properties of our cohomology groups. In particular, we construct a
long exact sequence which can be considered as a categorification
of the deletion-contraction rule of the chromatic polynomial. In
section 4, we show some computational examples. When the algebra
is $\Athree$, the cohomology groups can have a torsion of order
three. This is contrary to the cohomology groups in \cite{HR} and
in \cite{K}, where the computations suggest no odd torsion can
occur \cite{S}. We also show some computations when the algebra
has no grading, in which case $q\dim \A = \dim \A$ is an integer.
The Euler characteristic of the cohomology groups is therefore an
integer $P_G(\dim \A)$.  This should probably be compared to the
work of Eastwood and Huggett \cite{EH}, who constructed, for each
positive integer $\lambda $ and each graph $G$, a topological
space $M$ whose Euler characteristic is the integer
$P_G(\lambda)$. However, we don't know if there is any connections
between our work and theirs. We also make some comments on the
strength of the cohomology groups. In particular, we show that (a
twisted version) of our cohomology groups can be stronger than the
chromatic polynomial. The last section is an appendix in which we
classify rings whose additive group is $\mathbb{Z}\oplus
\mathbb{Z}$.

The authors wish to thank Mikhail Khovanov for his suggestions and comments.

\section{The Construction}
\subsection{The chromatic polynomial}
We recall some basic properties for the chromatic polynomial. Let
$G$ be a graph with vertex set $V(G)$ and edge set $E(G)$. For each
positive integer $\lambda$, let $\{ 1, 2, \cdots, \lambda\}$ be the
set of $\lambda$-colors. A \emph{$\lambda$-coloring} of $G$ is an
assignment of a $\lambda$-color to each vertex of $G$ such that
vertices that are connected by an edge in $G$ always have different
colors. Let $P_G(\lambda)$ be the number of $\lambda$-colorings of
$G$. It is well-known that  $P_G(\lambda)$ satisfies the
\emph{deletion-contraction relation}

\[ P_G(\lambda)=P_{G-e}(\lambda) - P_{G/e}(\lambda) \]
Furthermore, it is obvious that
\[ P_{N_n}(\lambda)=\lambda^n \mbox{ where $N_n$ is the graph with $n$
vertices and no edges.} \]

These two equations uniquely determines $P_G(\lambda)$. They also
imply that $P_G(\lambda)$ is always a polynomial of $\lambda$, known
as the {\it chromatic polynomial}.

There is another formula for $P_G(\lambda)$  that is useful for us. For each
$s\subseteq E(G)$, let $[G:s]$ be the graph whose vertex set is
$V(G)$ and whose edge set is $s$, let $k(s)$ be the number of
connected components of $[G:s]$. We have

\begin{equation} \label{State Sum 1}
 P_G(\lambda)=\mathrel{\mathop{\sum }\limits_{s\subseteq
E(G)}}(-1)^{\left| s\right| }\lambda ^{k(s)}.
\end{equation}

Equivalently, grouping the terms $s$ with the same number of edges
yields the following state sum formula

\begin{equation}
\label{State sum Chrom} P_G(\lambda)=\mathrel{\mathop{\sum
}\limits_{i\geq 0}}(-1)^{i}\mathrel{\mathop{\sum
}\limits_{s\subseteq E(G),\left| s\right| =i}}\lambda ^{k(s)}.
\end{equation}

\subsection{Graded algebras} We recall some definitions and specify
what kind of algebras we will work with.

\begin{definition}
Let $\R$ be a commutative ring with identity. An \emph{algebra} over
$\R$ is a ring $\A$  that is simultaneously a $\R$-module and such
that $r(ab)=(ra)b=a(rb)$ for all $r\in \R$ and $a, b\in \A$.
\end{definition}

\begin{definition}
A \emph{graded $\R$-algebra} $\A$ is an algebra with direct sum
decomposition $\A=\oplus_{i=0}^{\infty} A_i$ into $\R$-submodules
such that $a_i a_j \in A_{i+j}$ for all $a_i\in A_i$ and $a_j\in
A_j$. The elements of $A_j$ are called \emph{homogeneous} elements
of degree $j$.
\end{definition}

Most of our results can be stated for any integral domain $\R$
(see section 2.4). However, for simplicity, we will assume that
$R=\mathbb{Z}$, i.e. $\A$ is a ring itself. Furthermore, we want
each $A_i$ to have finite free rank, where the \emph{free rank} of
a $\Z$-module $M$ is defined to be $\dim_{\mathbb{Q}}(M
\otimes_{\Z} \mathbb{Q})$ and is denoted by $rank(M)$.

Thus, from now on, we will work with algebras that satisfy the
following conditions:

\begin{assumptions}
\label{assumptions}
 $\A=\oplus_{i=0}^{\infty} A_i$ is a commutative, graded
 $\Z$-algebra with $1$ such that each $A_i$ is free of finite rank.
\end{assumptions}

Note that these assumptions can sometimes be relaxed. For instance,
the construction can still be made even if there is no identity or
if the $A_i$'s are not free.

\bigskip

Since $\A$ can be consider as a $\mathbb{Z}$-module, the graded
dimension of $\A$ is defined using the definition below.

\begin{definition}
\label{graded dim Z-module} Let $M=\oplus_{j=0}^{\infty}M_{j}$ be
a graded $\mathbb{Z}$-module where $M_{j}$ denotes the set of
homogeneous elements of degree j of $M$. Assume that $rank (M_j)
<\infty$ for each $j$. The \emph{ graded dimension of }$M$ is the
power series
\[
q\dim M:=\sum _{j=0}^{\infty} q^{j}\: rank(M_{j}).
\]
\end{definition}

\subsection{The chain complex} Figure (\ref{The chain complex}) shows
what the chain complex will look like and the details can be found
after the figure. The diagram comes first because we thought it
might be helpful to have a picture of what is going on while reading
the formal definitions.

\begin{figure}[ht]
\begin{center}
\scalebox{.8}{\includegraphics{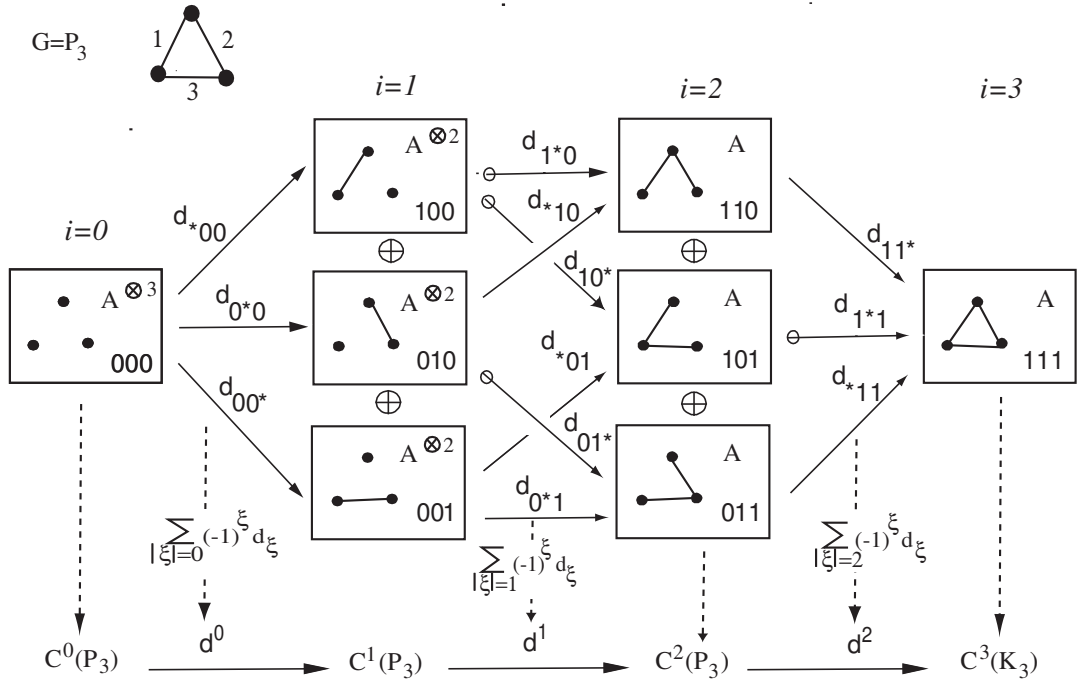}}
\caption{\textit{The chain complex when $G=P_3$}} \label{The chain
complex}
\end{center}
\end{figure}

 Let $G$ be a graph and $E=E(G)$ be the edge set of $G$. Let
$n=|E|$ be the cardinality of $E$. We fix an ordering on $E$ and
denote the edges by $e_1, \cdots, e_n$. Consider the
$n$-dimensional cube $\{ 0, 1\}^E = \{ 0, 1\}^n$. Each vertex
$\alpha$ of this cube corresponds to a subset $s=s_{\alpha}$ of
$E$, where $e_i \in s_{\alpha}$ if and only if $\alpha_i=1$. The
\emph{height} $|\alpha|$ of $\alpha$, is defined by $|\alpha|=\sum
\alpha_i$, which is also equal to the number of edges in
$s_{\alpha}$.

For each vertex $\alpha$ of the cube, we associate the graded
$\mathbb{Z}$-module $C^{\alpha}(G)$ as follows. Consider $[G:s]$,
the graph with vertex set $V(G)$ and edge set $s$. We assign a
copy of $\A$ to each component of $[G:s]$ and then taking tensor
product over the components. Let $C^{\alpha}(G)$ be the resulting
graded $\mathbb{Z}$-module, with the induced grading from $\A$.
Therefore, $C^{\alpha}(G) \cong \A^{\otimes k}$ where $k$ is the
number of components of $[G:s]$. We define the $i^{\mbox{\tiny
th}}$ chain group of our complex to be
$C^i(G):=\oplus_{|\alpha|=i} C^{\alpha}(G)$. Keep in mind that
$C^i(G)$ depends on the algebra ${\A}$. Thus one may want to
denote it by $C^i_{\A}(G)$. However, we will omit the letter
${\A}$ unless there is an ambiguity. Also, we sometimes
interchange the notions $\alpha$ and $s$. Thus $C^{\alpha}(G)$ is
sometimes denoted by $C^s(G)$. This should not cause any
confusion.

To define the differential maps $d^i$, we need to make use of the
edges of the cube $\{ 0, 1\}^E$. Each edge $\xi$ of $\{ 0, 1\}^E$
can be labeled by a sequence in $\{ 0, 1, *\}^E$ with exactly one
$*$. The tail of the edge is obtained by setting $*=0$ and the
head is obtained by setting $*=1$. The \emph{height} $|\xi|$ is
defined to bethe height of its tail, which  is also equal to the
number of 1's in $\xi$.

Given an edge $\xi$ of the cube, let $\alpha_1$ be its tail and
$\alpha_2$ be its head. The per-edge map $d_{\xi}:
C^{\alpha_1}(G)\rightarrow C^{\alpha_2}(G) $ is defined as
follows. For $j=1$ and 2, the $\Z$-module $C^{\alpha_j}(G)$ is
$\A^{\otimes k_j}$ where $k_j$ is the number of connected
components of $[G:s_j]$ (here $s_j$ stands for $s_{\alpha_j}$).
Let $e$ be the edge with $s_2=s_1\cup \{ e\}$.

If $e$ joins a component of $[G:s_1]$ to itself, then $k_1=k_2$
and the components of $[G:s_1]$ and the components of $[G:s_2]$
naturally correspond to each other. We let $d_{\xi}$ to be the
identity map.

 If $e$ joins two different components of $[G:s_1]$,
say $E_1$ to $E_2$ where $E_1, E_2, \cdots, E_{k_1}$ are the
components of $[G:s_1]$, then $k_2=k_1 - 1$ and the components of
$[G:s_2]$ are $E_1\cup E_2 \cup \{ e\}, E_3, \cdots, E_{k_1}$. We
define $d_{\xi}$ to be the identity map on the tensor factors
coming from $E_3, \cdots, E_{k_1}$, and $d_{\xi}$ on the remaining
tensor factors to be the multiplication map $\A \otimes \A
\rightarrow \A$ sending $x\otimes y$ to $ xy$.

Now, we define the differential $d^i:C^i(G) \rightarrow
C^{i+1}(G)$ by $d^i = \sum_{|\xi|=i} (-1)^{\xi} d_{\xi}$, where
$(-1)^{\xi}=-1$ (resp. 1) if the number of 1's in $\xi$ before $*$
is odd (resp. even).

We have

\begin{theorem}\label{chain complex}
 (a) $0 \rightarrow C^0(G)
\stackrel{d^0}{\rightarrow}C^1(G) \stackrel{d^1}{\rightarrow} \cdots
\stackrel{d^{n-1}}{\rightarrow} C^n(G) \rightarrow 0$ is a graded
chain complex whose differential
is degree preserving. \\
(b) The cohomology groups $H^i(G) (=H^i_{\cal A}(G))$ are
independent of the ordering of
the edges of $G$, and therefore are invariants of the graph $G$. \\
(c) The graded Euler characteristic of the chain complex is equal to
the chromatic polynomial of the graph $G$ evaluated at
$\lambda=q\dim \A$, i.e. $\chi _{q}(C)=\mathrel{\mathop{\sum
}\limits_{0\leq i\leq n}}(-1)^{i} q\dim
(H^{i})=\mathrel{\mathop{\sum }\limits_{0\leq i\leq n}} (-1)^{i}
q\dim (C^{i})=P_G(q\dim \A)$
\end{theorem}

\begin{proof} The proof is rather standard. We sketch the ideas here.

(a). To prove this defines a chain complex, we need to show that $d$
is a differential. That is, $d\circ d =0$. This is done in two
steps. First, we verify that the maps $d_{\xi}$ makes the cube
commutative, a fact follows from the associativity of the algebra.
Second, the signs $(-1)^{\xi}$ in $d$ allow us to cancel out all
terms in $d\circ d$. Thus $d\circ d=0$.

 To show that $d$ is degree
preserving, we note that the multiplication map on \A\  is always
degree preserving, which then implies each map $d_{\xi}$ is degree
preserving, and therefore so is $d$.

(b). The proof is similar to the one for $\Atwo$ in \cite{HR}.
Each permutation of the edges of $G$ is a product of
transpositions of the form $(k, k+1)$. An explicit isomorphism can
be constructed for each such transposition. In fact, this shows
that the isomorphism class of the chain complex is an invariant of
the graph.

(c). First, a standard homological algebra argument
 shows that $\sum _{0\leq i\leq n}(-1)^{i} q\dim
(H^{i})=\sum_{0\leq i\leq n} (-1)^{i} q\dim (C^{i})$.  Next, we
note that $q \dim C^{\alpha}(G) =q\dim \A^{\otimes k}=(q\dim
\A)^k$ where $k$ is the number of connected components of $[G:s]$.
Taking direct sum over $s\subseteq E(G), |s|=i$, and then taking
alternating sum over $i$, we obtain the equation $\sum_{0\leq
i\leq n} (-1)^{i} q\dim (C^{i})= P_G(q\dim \A)$.
\end{proof}

\begin{remark}
(a) The above graded chain complex can be turned into a bi-graded
chain complex. Let $C^{i,j}(G)$ be the subgroup of $C^i(G)$
consisting of homogeneous elements with degree $j$. For each $j$ we
have a chain complex

\[ 0 \rightarrow C^{0,j}(G)
\stackrel{d^{0,j}}{\rightarrow}C^{1,j}(G)
\stackrel{d^{1,j}}{\rightarrow} \cdots
\stackrel{d^{n-1,j}}{\rightarrow} C^{n,j}(G) \rightarrow 0 \]

The direct sum of these chain complexes, with the obvious gradings,
is equal to the chain complex in Theorem \ref{chain complex}.

(b) Our chain complexes can also be described in terms of enhanced
states. One defines an \emph{enhanced state} $S$ of $G$  to be a
pair $(s,c)$ where $s\subseteq E(G)$ and $c$ is an assignment of an
element of $\A$ to each connected component of $[G:s]$. Identify $S$
with the element $c(E_1) \otimes \cdots \otimes c(E_k)$ of
$C^{s}(G)=\A^{\otimes k}$, where $E_1, \cdots, E_k$ are components
of $[G:s]$. Thus $C^i(G)$ is generated by states with $|s|=i$. When
each $c(E_i)$ is a homogeneous element of $\A$, we say the coloring
$c$ and the enhanced state $S$ are homogeneous, and we define its
degree to be $j(S)=\sum_i \deg c(E_i)$. It is easy to see that
$C^{i,j}$ above is generated by all homogeneous enhanced states $S$
with $i(S)=i, j(S)=j$. The differential of each enhanced state is
then defined to be the operation of adding each edge not in $s$,
adjusting the coloring $c$ using the multiplication on $\A$ or the
identity map, and then taking the summation over the edges in
$E(G)-s$ with appropriate $\pm 1$ signs in front of each term.

\end{remark}

\subsection{Some variations in the construction}
Some variations can be introduced in our construction.

\bigskip

(a) As we discussed in 2.2, our coefficient ring $\mathbb{Z}$ can
be replaced by any integral domain $R$, i.e. a commutative ring
with $1\neq 0$ and no zero divisors. By analogy with the abelian
group case, we define the \emph{free rank} of a module $M$ over an
integral domain $\R$ to be $rank(M):=dim_Q(M \otimes_R Q)$
where $Q$ is the field of fractions of $R$. 
 Note
that $M \otimes_R Q$ is
a vector space over $Q$ so its dimension is well defined.

Let $A=\oplus_{i=0}^{\infty} A_i$ be a graded algebra over an
integral domain $R$ such that each $A_i$ has finite free rank.
Using $A$ as the building block for our chain $R$-modules (the
counterpart of the chain groups in the $R$-module case), the same
construction can still be made and $\sum_{0\leq i\leq n} (-1)^{i}
q\dim (C^{i})=P_G(q\dim \A)$ where the graded dimension of $A$ is
obviously defined to be the power series $q\dim \A =\sum_{i}q^i
rank A_i$.

The proof that $\sum_{0\leq i\leq n} (-1)^{i} q\dim
(H^{i})=\sum_{0\leq i\leq n} (-1)^{i}q\dim (C^{i})$
 is the
same as in the abelian case, we just need to prove the counterpart
in the $R$-module category of a well known result for abelian
groups. Namely, it is enough to show that if $0\rightarrow L
\rightarrow M \rightarrow N \rightarrow 0$ is an exact sequence of
$R$-modules, then $rank M= rank L + rank N$. This is true because
the quotient field of an integral domain $R$ is a flat $R$-module
(\cite{R}, p86). Hence taking the tensor product of the exact
sequence $0\rightarrow L \rightarrow M \rightarrow N \rightarrow
0$ by $Q$ yields an exact sequence of $R$-modules that turns out
to be an exact sequence of $Q$-vector spaces and the result
follows.

Some properties, e.g. the result (\ref{summand}) on pendant edge
in the next section, can be extended to the case when $\R$ is a
principal idea domain.

\bigskip

(b) We can define a twisted version of the differential map $d$ as
follows. Let $f: \A \rightarrow \A$ be a degree preserving
homomorphism of the algebra $\A$. That is, $f$ is linear,
$f(yz)=f(y)f(z)$ for all $y,z$, and $deg (f(y))=deg (y)$ for all
homogeneous elements $y$. Let $\xi$ be an edge of the cube we had
before. Let $\alpha_1$ be its tail, $\alpha_2$ be its head, and
$e$ be the edge of $G$ corresponding to $\xi$ so that $s_2=s_1\cup
\{ e\}$. As before, the $\Z$-module $C^{\alpha_j}(G) = \A^{\otimes
k_j}$ where $j=1, 2$ and $k_j$ is the number of connected
components of $[G:s_j]$ (here $s_j=s_{\alpha_j}$).

If $e$ joins a component of $[G:s_1]$ to itself, then $k_1=k_2$
and the components of $[G:s_1]$ and the components of $[G:s_2]$
naturally correspond to each other. We let $d_{\xi}=f^{\otimes
k_1}: \A^{\otimes k_1}\rightarrow \A^{\otimes k_2}$ (note that
$k_1=k_2$).

If $e$ joins two different components of $[G:s_1]$, we define
$d_{\xi}$ the same way as before. This yields the differential
$d^i:C^i(G) \rightarrow C^{i+1}(G)$ using the same formula: $d^i =
\sum_{|\xi|=i} (-1)^{\xi} d_{\xi}$.

If $f$ is the identity map, this is the same as our construction
above. Otherwise, the constructions can be different. For example,
if $f=0$ and $G$ is the graph with one vertex and $n$ loops, then
the differential $d=0$. Thus $H^i(G)\cong C^i(G)$ which are
nontrivial for $i=0, 1, \cdots, n$. On the other hand, the
(untwisted) cohomology groups all vanishes when $G$ has a loop
(see Corollary \ref{loops} in the next section).

\section{Some Properties}

In this section, we discuss some basic properties of the cohomology
groups. Most of these are parallel to the case when the algebra is
$\Ai$ whose discussion was given in \cite{HR}. The most interesting
property is a long exact sequence, which can be considered as a
categorification for the deletion-contraction rule for the chromatic
polynomial. \\

\subsection{An exact sequence} The long exact sequence comes from a
short exact sequence of graded chain homomorphisms

\[0\rightarrow C^{i-1}(G/e)\overset{\alpha}{\rightarrow} C^{i}(G)
\overset{\beta}{\rightarrow} C^{i}(G-e) \rightarrow 0 \] which we
explain here. Basically $\alpha$ is the map that recovers the edge
$e$, and $\beta$ is the projection map that kills every state
containing $e$. A more precise description is given below. First,
we order the edges of $G$ so that $e$ is the last edge. This
induces natural orderings on the edge sets of $G/e$ and $G-e$ by
deleting $e$ from the list. For each $s\subseteq E(G/e)$, let
$\tilde{s}=s\cup \{ e\}$. Then $\tilde{s}\subseteq E(G)$. Recall
that $C^{s}(G/e)$ (resp. $C^{\tilde{s}}(G)$) is the tensor product
of $\A$ taken over components of $[G/e:s]$ (resp.
$[G:\tilde{s}]$). The components of $[G/e:s]$ and the components
of $[G:\tilde{s}]$ are the same except for the one involving $e$
in which case they are related by a contraction of $e$.  Thus we
have $C^s(G/e)\cong C^{\tilde{s}}(G)$ via a natural isomorphism,
since the tensor factors naturally correspond to each other. Let
$\alpha|_{C^s(G/e)}: C^s(G/e)\rightarrow C^{\tilde{s}}(G)$ be this
isomorphism. Taking direct sum over $s$, we obtain the
homomorphism $\alpha: C^{i-1}(G/e) \rightarrow C^i(G)$.
 Next, we explain the
map $\beta: C^{i}(G) \rightarrow C^{i}(G-e)$. We have $C^i(G)=
\oplus_{s\in E(G),|s|=i} C^s(G)$. If $e\not \in s$, $s$ is
automatically a subset of $E(G-e)$. We have $C^s(G)=C^{s}(G-e)$
since the graphs $[G:s]$ and $[G-e:s]$ are identical. The map $\beta
$ is the identity map from $C^s(G)$ to $C^s(G-e)$. If $e\in s$, we
let $\beta|_{C^s(G)}$ be the zero map. Taking direct sum over $s$
with $|s|=i$, we obtain the map $\beta: C^i(G)\rightarrow C^i(G-e)$.
A standard diagram chasing argument shows that this defines a short
exact sequence of chain complexes. Thus we have

\begin{theorem}
Let $G$ be a graph, and $e$ be an edge of $G$. \\
 (a) For each $i$, there is a short exact sequence of graded chain homomorphisms:
$0\rightarrow C^{i-1}(G/e)\overset{\alpha}{\rightarrow} C^{i}(G)
\overset{\beta}{\rightarrow} C^{i}(G-e) \rightarrow 0$, and
therefore by the zig-zag lemma, \\
(b) it induces a long exact sequence of cohomology groups:
$0\rightarrow H^{0}(G)\overset{\beta ^{\ast }}{\rightarrow
}H^{0}(G-e) \overset{\gamma ^{\ast }}{\rightarrow
}H^{0}(G/e)\overset{\alpha ^{\ast }}{ \rightarrow
}H^{1}(G)\overset{\beta ^{\ast }}{\rightarrow} \ldots
\rightarrow H^{i}(G)\overset{\beta ^{\ast }}{\rightarrow }H^{i}(G-e)%
\overset{\gamma ^{\ast }}{\rightarrow }H^i(G/e)\overset{\alpha
^{\ast }}{ \rightarrow
} \ldots \ $ \\
\end{theorem}

Taking the alternating sum of the graded dimensions in the above
long exact sequence, we obtain the deletion-contraction rule. It is
in this sense that the long exact sequence is considered as a
categorification of the deletion-contraction rule.

\bigskip

It is useful to understand the following geometric description of
the maps $\alpha^{\ast}, \beta^{\ast},$ and $ \gamma^{\ast}$:
$\alpha ^{\ast }$ expands the edge $e$ and keeps in the same
coloring, $\beta ^{\ast }$ is the projection map that kills every
state containing $e$. For $\gamma^*$, we add the edge $e$, contract
it to a point, assign a natural coloring described below and then
multiply by $(-1)^i$. The corresponding coloring is the same if $e$
connects a component to itself. If $e$ connects two separate
components, say $E_1$ and $E_2$, then new coloring on the component
$E_1\cup E_2\cup \{ e\}$ is the product $c(E_1)c(E_2)$ and is the
same as $c(E_i)$ for all other components $E_i$.

Other basic properties can follow either from this exact sequence or
from the definition. For example, we have

\subsection{Graphs with loops or multiple edges}

\begin{corollary}\label{loops}
(a) \label{graphs with loops} If a graph has a loop then all the
cohomology groups are trivial.\\
(b) \label{graphs with multiple edges} The cohomology groups are
unchanged if you replace all the multiple edges of a graph by single
edges.
\end{corollary}

\begin{proof} (a). In the long exact sequence

\[ \ldots {\rightarrow }H^{i-1}(G-e)%
\overset{\gamma ^{\ast }}{\rightarrow }H^{i-1}(G/e)\overset{\alpha
^{\ast }}{ \rightarrow} H^{i}(G)\overset{\beta ^{\ast }}{\rightarrow }H^{i}(G-e)%
\overset{\gamma ^{\ast }}{\rightarrow }H^i(G/e)\overset{\alpha
^{\ast }}{ \rightarrow } \ldots \ \]

\noindent we have $G/e = G-e$ and the map $\gamma^{\ast}$ is the
identity map multiplied by $(-1)^i$. It follows that $H^i(G)=0$ for
each $i$.

(b) Let $e_1$ and $e_2$ be two edges connecting the same pair of vertices in $G$.
 In the exact sequence
\[ \rightarrow H^{i-1}(G/e_2)\overset{\alpha^{\ast }}{ \rightarrow}
H^{i}(G)\overset{\beta ^{\ast }}{\rightarrow }H^{i}(G-e_2)%
\overset{\gamma ^{\ast }}{\rightarrow }H^i(G/e_2)\rightarrow \] the
graph $G/e_2$ contains a loop coming from $e_1$. Therefore
$H^{i-1}(G/e_2)=H^i(G/e_2)=0$. It follows that $H^{i}(G)\cong
H^{i}(G-e_2)$. One can repeat this process until all redundant edges
in $G$ are removed.
\end{proof}

\subsection{Adding a pendant edge.}  Next, we consider the effect of
adding a pendant edge to a graph. Recall that a {\em pendant vertex}
in a graph is a vertex of degree one, and a {\em pendant edge} is an
edge connecting a pendant vertex to another vertex. Let $e$ be a
pendant edge in a graph $G$, then
$P_G(\lambda)=(\lambda-1)P_{G/e}(\lambda)$. An analogous equation on
the cohomology level is given below, as long as our algebra $\A$ has
an identity. First we prove an algebraic lemma.


\begin{lemma}\label{summand}
Let $\A$ be an algebra over $\mathbb{Z}$ that with identity
$1_{\A}$. Assume that the additive group $(\A , +)$ is a free
abelian group of finite rank. Then $1_{\A }$ generates a direct
summand of the abelian group $(\A, +)$.
\end{lemma}

\begin{proof} Let $d$ be a positive integer such that
$1_{\A}=de_1$ and $\{ e_1, \cdots, e_k\}$ is a basis of the
abelian group $(\A, +)$ (see \cite{H} p73). It enough to show that
$d=1$. We have $1_{\A}=de_1=d(1_{\A}e_1)=d^2 (e_1^2)$. We write
$e_1^2$ as a linear combination of the basis elements:
$e_1^2=b_1e_1+\cdots + b_ke_k$. This implies $1_{\A}=d^2b_1 e_1 +
\cdots d^2 b_k e_k$. On the other hand, we have $1_{\A}=d e_1$.
The uniqueness of coefficients implies that $d^2b_1=d$, which
implies $d=1$.
\end{proof}

\bigskip
By the lemma, if our graded algebra $\A$ has an identity, then $\A
= \mathbb{Z}1_{\A\ } \oplus \A'$ as a $\mathbb{Z}$-module, where
$\mathbb{Z}1_{\A }$ is generated by the identity of $\A$ and $\A'$
is a submodule of $\A$. We simply apply the lemma to our graded
algebra $\A$. If $\A$ has infinite dimension, we apply the lemma
to the subalgebra $A_0$ instead which is free of finite rank by
assumptions (\ref{assumptions}).

\begin{proposition}\label{pendant edge} Let $\A$ and $\A'$ be as above.
Then $H^i(G)\cong H^i(G/e)\otimes \A'$ where $e$ is a pendant edge
of $G$.
\end{proposition}

\begin{proof}
Consider the operations of contracting and deleting $e$ in $G$.
Denote the graph $G/e$ by $G_1$. We have $G/e=G_1,$ and
$G-e=G_1\sqcup \{v\}$, where $v$ is the end point of $e$ with $\deg v = 1$.
Consider the exact sequence

\begin{equation*}
\cdots \rightarrow H^{i-1}(G_1\sqcup
\{v\})\overset{\gamma^*}{\rightarrow } H^{i-1}(G_1)
 \overset{\alpha^*}{\rightarrow } H^{i}(G)\overset{\beta^*}{\rightarrow
} H^{i}(G_1 \sqcup \{v\})\overset{\gamma^*}{\rightarrow }
H^{i}(G_1)\rightarrow \cdots
\end{equation*}

We need to understand the map
\begin{equation*}
H^{i}(G_1\sqcup \{v\})\overset{\gamma^*}{\rightarrow }H^{i}(G_1)
\end{equation*}

It is easy to compute the cohomology groups for the one point
graph $\{ v\}$(see the first example in the next section). We have
$H^0(\{ v\})\cong \A$ and $H^i(\{ v\})=0$ for all $i>0$. Thus, the
K\"{u}nneth type formula (\ref{disjoint union}) below implies

\begin{equation*}
H^{i}(G_{1}\sqcup \{v\})\cong H^{i}(G_{1})\otimes \A
\end{equation*}

\noindent by a natural isomorphism $h_*$. Since, $\A
=\mathbb{Z}1_{\A} \oplus \A'$, we identify $H^{i}(G_1\sqcup \{v\})$
with $H^{i}(G_1) \otimes (\Z \oplus \A')$. The map $\gamma^{\ast} :
H^{i}(G_1\sqcup \{v\})\rightarrow H^{i}(G_1)$ sends $x\otimes 1$ to
$(-1)^i x$. In particular, $\gamma^*$ is onto. Therefore, the above
long exact sequence becomes a collection of short exact sequences

\begin{equation}
0\rightarrow H^{i}(G)\overset{\beta^*}{\rightarrow }
H^{i}(G_{1}\sqcup \{v\})\overset{\gamma^*}{\rightarrow }
H^{i}(G_{1})\rightarrow 0
\end{equation}

Hence, $H^i(G)\cong \ker \gamma^* $. We define a homomorphism:

\[ f: H^i(G_1)\otimes \A' \rightarrow \ker \gamma^* \mbox {  by } \\
 f(x\otimes a')=x\otimes a' -(-1)^i \gamma^*(x\otimes a') \otimes 1 \]

One checks that $f$ is an isomorphism of $\Z$-modules. Therefore,
$H^i(G)\cong \ker \gamma^* \cong H^i(G_1)\otimes \A'$.
\end{proof}

\begin{remark}\label{pendant edge rmk}
The isomorphism $f: H^i(G_1)\otimes \A' \rightarrow H^i(G)$
defined above can be visualized as follows. Each cycle in
$C^i(G_1)$ is a linear combination of terms of the form
$\vcenter{\hbox{\epsfig{file=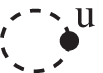}}}$, where the vertex to which
$e$ collapses is labelled with $u\in \A$. Such a term becomes
$\vcenter{\hbox{\epsfig{file=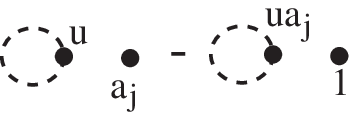}}}$ under the action of
$f$, where the element of $\A$ written close to a vertex indicates
the label that has been assigned to this component.

\end{remark}

\subsection{Disjoint union}

Finally, we state a K\"{u}nneth theorem type formula for our
cohomology groups under disjoint union. It will be used in our
computation in the next section.

\begin{proposition}
\label{disjoint union}Let $\A$ be an algebra satisfying the
assumptions (\ref{assumptions}). For each $i\in \mathbb{N}$, we
have:
\begin{equation*}
H_{\A}^{i}(G_{1}\sqcup G_{2})\cong \left[ \mathrel{\mathop{\oplus
}\limits_{p+q=i}} H_{\A}^{p}(G_{1})\otimes
H_{\A}^{q}(G_{2})\right] \oplus \left[ \mathrel{\mathop{ \oplus
}\limits_{p+q=i+1}} H_{\A}^{p}(G_{1})*H_{\A}^{q}(G_{2})
 \right]
 \end{equation*}
 where * denotes the torsion product of
two abelian groups.
\end{proposition}

\begin{proof}
This is a corollary of K\"{u}nneth theorem \cite{Munkres}, since the
chains complexes $ \mathcal{C}(G_{1})$\ and $\mathcal{C}(G_{2})$\
are free.
\end{proof}

\section{Some Computations}



We show some computational examples. The examples for null graphs
and trees follow easily from the previous section. The results are
needed for other computations. The examples for the graph $P_3$,
``the triangle", with different algebras show that there can be a
variety of torsions. Our final example shows the (twisted)
cohomology theory (see section 2.4) can be stronger than the
chromatic polynomial, by varying the algebra $\A$ and the
homomorphism $f$.

\begin{example}\label{N_n}
 Let  $G=N_n$ be the order $n$ null graph. That is, the graph with
$n$ vertices and no edges. Then $C^0(G)\cong \A^{\otimes n}$ and
$C^i(G)=0 $ for $i>0$. It follows that $H^0(G)\cong \A^{\otimes n}$
and $H^i(G)=0$ for all $i>0$.
\end{example}

\begin{example}\label{union}
Given any graph $G$, let $G\sqcup N_n$ be the graph obtained by
adding $n$ isolated vertices to $G$. By Proposition \ref{disjoint
union} and Example \ref{N_n},  $H^i(G\sqcup N_n)\cong H^i(G) \otimes
\A^{\otimes n}$ for all $i\geq 0$.
\end{example}

\begin{example}\label{T_n}
Let $G=T_n$ be a tree with $n+1$ vertices and $n$ edges. We assume
$\A$ has an identity and therefore $\A \cong \mathbb{Z}\oplus \A'$
for some $\mathbb{Z}$-submodule $\A'$. By Proposition \ref{pendant
edge}, $H^0(T_n)\cong \A \otimes \A'^{\otimes n}$ and $H^i(T_n)=0$
for all $i>0$.
\end{example}



Recall that $\cdot \{\ell\}$ is the \emph{degree shift} operation on
graded $\Z$-modules that increases the degree of all homogeneous
elements by $\ell$. For instance, in the example below,
$\mathbb{Z}^3\{ 4\}$ denotes the graded $\mathbb{Z}$-module
 $\mathbb{Z}\oplus \mathbb{Z}\oplus \mathbb{Z}$ where the three generators have degree 4.

\begin{example}
Let $G=P_3$, the polygon graph on 3 vertices (see Figure (\ref{The
chain complex})). Let
 $\A\ = \Athree$, we will show below that
\[ H^0(P_3) \cong \mathbb{Z} \{3\} \oplus \mathbb{Z}^3 \{4\} \oplus
\mathbb{Z}^3 \{5\} \oplus \mathbb{Z}\{6\} \]

\[ H^1(G)\cong \mathbb{Z}\{1\} \oplus \mathbb{Z}\{2\} \oplus
\mathbb{Z}_3 \{3\} \] and $H^i(G)=0$ for $i>1$.

\end{example}

The computation is done using the exact sequence. We start with
$N_1$, and then add edges until we have $P_3$.

If $G=N_1$, $H^0(G)\cong \A$ with generators being $\{ x^r | r=0, 1,
2 \}$.

If $G=T_1$, the tree with one edge, Proposition \ref{pendant edge}
and Remark \ref{pendant edge rmk} imply that $H^0(G)\cong \A\
\otimes \A'$ with generators $\{ e_{rs} | r=0,1,2, s=1,2\}$, where
$e_{rs}=x^r \otimes x^s - x^{r+s}\otimes 1$ in $H^0(G)$. For
simplicity, we drop the tensor notation and denote $e_{rs}$ by
$e_{rs}=x^rx^s-x^{r+s}1$. A picture of $e_{rs}$ is $e_{rs}=
\vcenter{\hbox{\epsfig{file=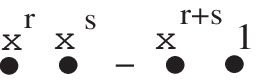}}}$

If $G=T_2$, the tree with two edges, the same argument shows that a
set of generators of $H^0(G)\cong H^0(T_1)\otimes \A'$ is $\{
e_{rst} | r=0, 1, 2, s=1,2$ and $t=1,2\}$  where
$e_{rst}=(x^rx^sx^t-x^rx^{s+t}1)-(x^{r+s}1x^t-x^{r+s}x^t1)=x^rx^sx^t-x^rx^{s+t}1-x^{r+s}1x^t+x^{r+s}x^t1$.
A picture of $e_{rst}$ is
$e_{rst}=\vcenter{\hbox{\epsfig{file=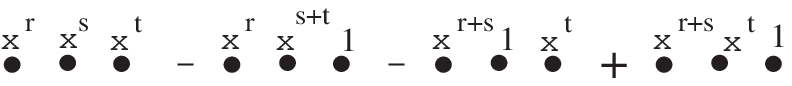}}}$.

Now, let $G=P_3$. Let $e$ be an edge of $P_3$. The exact sequence on
$(P_3, e)$ gives

\[ 0\rightarrow H^0(P_3) \overset{\beta^*}{\rightarrow} H^0(P_3-e)
\overset{\gamma^*}{\rightarrow} H^0(P_3/e)
\overset{\alpha^*}{\rightarrow} H^1(P_3) \rightarrow 0 \]

\noindent
where $H^0(P_3-e)=H^0(T_2)$ is freely generated by \\
$\{ e_{011}, e_{012}, e_{021}, e_{111}, e_{022}, e_{112}, e_{121}, e_{211}, e_{122}, e_{212}, e_{221}, e_{222} \}$,\\
and $H^0(P_3/e)\cong H^0(T_1)$ is freely generated by $\{e_{01},
e_{02}, e_{11}, e_{12}, e_{21}, e_{22}\}$. The map $\gamma^*$ sends
$a\otimes b \otimes c$ to $ac\otimes b$. Thus
$\gamma^*(e_{rst})=x^{r+t}x^s-x^rx^{s+t}-x^{r+s+t}1+x^{r+s}x^t
=(x^{r+t}x^s-x^{r+s+t}1)-(x^rx^{s+t}-x^{r+s+t}1)+(x^{r+s}x^t-x^{r+s+t}1)
=e_{r+t,s}-e_{r,s+t}+e_{r+s,t}$. This gives \\

$\indent
\gamma^*(e_{011})=e_{11}-e_{02}+e_{11}=2e_{11}-e_{02} \\
\indent \left\{
\begin{array}{ll} \gamma^*(e_{012})=e_{21}-e_{03}+e_{12}=e_{12}+e_{21}
\mbox{ (here $e_{03}=0$ since $x^3=0$, same for $e_{13}$ etc)} &  \\
 \gamma^*(e_{021})=e_{12}-e_{03}+e_{21}=e_{12}+e_{21} & \\
  \gamma^*(e_{111})=e_{21}-e_{12}+e_{21}=-e_{12}+2e_{21} &
\end{array} \right. \\
 \indent \left\{ \begin{array}{l}
\gamma^*(e_{022})=e_{22}-e_{04}+e_{22}=2e_{22} \\
\gamma^*(e_{112})=e_{31}-e_{13}+e_{22}=e_{22} \\
\gamma^*(e_{121})=e_{22}-e_{13}+e_{31}=e_{22}\\
\gamma^*(e_{211})=e_{31}-e_{22}+e_{31}=-e_{22}\\ \end{array} \right.
\\
\indent \left\{ \begin{array}{l} \gamma^*(e_{122})=e_{32}-e_{14}+e_{32}=0 \\
\gamma^*(e_{212})=e_{41}-e_{23}+e_{32}=0\\
\gamma^*(e_{221})=e_{32}-e_{23}+e_{41}=0\end{array} \right.
\\ \indent \indent \gamma^*(e_{222})=e_{42}-e_{24}+e_{42}=0 $ \\

Here, we group the basis elements into five groups according to
their degrees (i.e. $r+s+t$ in the domain). This breaks the
$12\times 6$ presentation matrix into five blocks of submatrices.
From there, we obtain the cohomology groups of $P_3$ using
 $H^0(P_3)\cong \ker \gamma^*$ and $H^1(P_3)\cong
H^0(P_3/e)/\mbox{Im}\gamma^*$.


\bigskip

In the next two examples, let us consider the special case when
$\A\ = A_0$, i.e., every element in \A\ has degree 0. The graded
dimension of \A\ is then an integer, namely $\dim \A$. Thus the
Euler characteristic of our cohomology groups is the integer
$P_G(\lambda )$ where $\lambda = \dim \A $. We study two examples,
the first one is when $\dim \A\ =1$, the second is when $\dim \A\
= 2$.


\begin{example}
Let $\A=\Z$ with the usual ring structure. We have $q\dim \A =1$,
$P(G,1)=1$ if $G$ has no edge and 0 otherwise. For the cohomology
groups, we have
\[
H^i(G)\cong \left\{
\begin{array}{ll} \mathbb{Z} & \mbox{ if  $i=0$ and $G$ has no edge, }\\
0 & \text{ otherwise.}
\end{array} \right.
\]

This can be easily proved by inducting on the number of edges and using the long exact sequence.
\end{example}

Next, we consider the case when $\dim \A = 2$. We assume that $\A$
has an identity. By Lemma \ref{summand}, the abelian group $(\A, +)$
is generated by $1$ and $x$ where $1$ is the identity of $\A$. Let
us consider various ring structures on $\A$. The multiplication $*$
satisfies $1*1=1, 1*x=x*1=x$, and $x*x=a1+bx$ where $a, b$ are two
fixed integers. In the appendix, we will show that the isomorphism
type of such a ring depends on $b\pmod{2}$ and $b^2+4a$.

\begin{example}\label{ring_exa}
Let \A\ be the ring above. Then

$H^0(P_3)\cong \left\{ \begin{array}{ll} \Z  & \mbox{ if $b^2+4a=0$,} \\
0 & \mbox{ otherwise.} \end{array} \right.$

$H^1(P_3)\cong \left\{ \begin{array}{ll} \Z_2\oplus \Z  & \mbox{ if $b^2+4a=0$,} \\
\Z_{|b^2+4a|} & \mbox{ if $b^2+4a \neq 0$ and $b$ is odd} \\
\Z_2\oplus Z_{\frac{|b^2+4a|}{2}} & \mbox{ if $b^2+4a\neq 0$ and $b$
is even}
\end{array} \right. $

$H^i(P_3)=0$ for all $i>1$.

\end{example}
This example shows that all torsions that are not of the form $3+4k$
can occur.

\medskip

Again, the computation is based on the exact sequence.

If $G=N_1$, the graph with one vertex and no edge, then $H^0(G)\cong \A$ with
generators being $1$ and $x$.

If $G=T_1$, the tree with one edge, Proposition \ref{pendant edge}
implies that $H^0(G)\cong \A\  \otimes Zx\cong \A$. By Remark
\ref{pendant edge rmk}, $H^0(G)$ has the following basis:
$e_1=1\otimes x-x\otimes 1=1x-x1, e_2=x\otimes x-x*x\otimes
1=x\otimes x-(a1+bx)\otimes 1 =x\otimes x-a1\otimes 1 -bx\otimes
1=xx-a11-bx1$. Here, for simplicity of notation, we drop the
tensor product symbol $\otimes$ again. Thus $1x$ stands for
$1\otimes x$.

If $G=T_2$, the tree with two edges, the same argument shows that
$H^0(T_2)\cong \A$. To describe the basis, we denote the three
vertices of $T_2$ by 1, 2, and 3 with 2 being the middle vertex with
degree two.
 A basis of $H^0(T_2)$ is:\\
$f_1=(1\otimes x \otimes x - 1\otimes x*x\otimes 1) - (x\otimes 1\otimes x-x\otimes 1*x\otimes 1)\\
=(1xx)-1(a1+bx)1-(x1x)+(xx1)=(1xx)-a(111)-b(1x1)-(x1x)+(xx1)$,\\
$f_2=(xxx)-[x(x*x)1]-a(11x-1x1)-b(x1x-xx1)\\=(xxx)-x(a1+bx)1-a(11x)+a(1x1)-b(x1x)+b(xx1)\\
=(xxx)-a(x11)-b(xx1)-a(11x)+a(1x1)-b(x1x)+b(xx1)$.

Now, let $G=P_3$. Let $e$ be an edge of $P_3$. The exact sequence on
$(P_3, e)$ gives

\[ 0\rightarrow H^0(P_3) \overset{\beta^*}{\rightarrow} H^0(P_3-e) \overset{\gamma^*}{\rightarrow} H^0(P_3/e)
\overset{\alpha^*}{\rightarrow} H^1(P_3) \rightarrow 0 \]

\noindent we have $H^0(P_3-e)=H^0(T_2)\cong \Z \oplus \Z$ with $\{
f_1, f_2\}$ being a basis, and $H^0(P_3/e)\cong H^0(T_1)\cong \Z
\oplus \Z$ with $\{ e_1, e_2\}$ being a basis. The map $\gamma^*$
sends $a\otimes b\otimes c$ to  $ac\otimes b$. Thus \\
$ \gamma^*(f_1)= \gamma^*((1xx)-a(111)-b(1x1)-(x1x)+(xx1))\\
=xx-a11-b1x-(x^2)1+xx =2(xx)-a(11)-b(1x)-(a1+bx)1
\\ =2[(xx)-a(11)-b(x1)]-[b(1x)-b(x1)]\\
=-be_1+2e_2.\\
\gamma^*(f_2)=\gamma^*((xxx)-a(x11)-b(xx1)-a(11x)+a(1x1)-b(x1x)+b(xx1)) \\
=(x*x)x-a(x1)-b(xx)-a(x1)+a(1x)-b((x*x)1)+b(xx)\\
=(a1+bx)x-2a(x1)-b(xx)+a(1x)-b((a1+bx)1)+b(xx)\\
=a(1x)+b(xx)-2a(x1)-b(xx)+a(1x)-ab(11)-b^2(x1)+b(xx)\\
=b(xx)-2a(x1)-b^2(x1)+2a(1x)-ab(11)\\
=2a(1x-x1)+b(xx)-ab(11)-b^2(x1)\\
=2a e_1 + b e_2 $

Therefore, the linear map $\gamma^*$ is given by the matrix $\left(
\begin{array}{ll} -b & 2 \\ 2a & b \end{array} \right)$.
Standard linear algebra then implies the computation results.




\bigskip

In our final example, we consider the twisted cohomology groups
described in section 2.4, and show that they can be stronger than
the chromatic polynomial. Of course, as we noted earlier in section
2.4, graphs with one vertex and $n$ loops provide such examples.
They all have zero chromatic polynomials. The example below contains
no loops, and their chromatic polynomials are nonzero.

\begin{example}
 Let $G_1$ be two copies of $P_3$ glued together at one vertex,
 $G_2$ be two copies of $P_3$ glued together at one edge, plus a pendant edge
 (here $P_3$ is the polygon graph on 3 vertices, see Figure (\ref{The chain complex}).
 It is easy to check that $G_1$ and $G_2$ share the same chromatic polynomial:

 \[ P_{G_1}=P_{G_2}=\lambda(\lambda-1)^2(\lambda-2)^2 \]

On the other hand, they often have non-isomorphic chain complexes,
and sometimes non-isomorphic twisted cohomology groups. For a
concrete example, let our algebra $\A$ be $\mathbb{Z}x$ with $\deg
x=1$ and $x*x=0$. Let $f=0$ be the homomorphism from $\A$ to $\A$.
Then the differential $d=0$. Therefore $H^i(G)\cong C^i(G)$ for all
$i$ and $G=G_1, G=G_2$. The chain groups $C^i$ are computed below.
We have

For $G_1$, $C^0\cong \mathbb{Z}\{ 5 \}, C^1\cong \mathbb{Z}^6\{ 4\},
C^2\cong \mathbb{Z}^{15}\{ 3\}, P_3\cong \mathbb{Z}^2\{ 3\}\oplus
\mathbb{Z}^{18}\{ 2\}, C^4\cong \mathbb{Z}^7\{ 2\}\oplus
\mathbb{Z}^8\{ 1\},
C^5=\mathbb{Z}\{ 2\}\oplus \mathbb{Z}^5\{ 1\}, C^6=Z\{ 1\}$,\\
 For $G_2$, $C^i$ are the same except $C^4\cong \mathbb{Z}^6\{ 2\}\oplus \mathbb{Z}^9\{ 1\},
 C^5=\mathbb{Z}^6\{ 1\}$. Here $\mathbb{Z}^2\{ 3\}$ denotes the graded $\mathbb{Z}$-module
 $\mathbb{Z}\oplus \mathbb{Z}$ with both generators having degree 3.

This shows that $C^i(G_1)\not \cong C^i(G_2)$ for $i=4, 5$, and
therefore their cohomology groups $H^i$ are different for $i=4$ and
5.
\end{example}

\section{Appendix - Ring structures on $\mathbb{Z}\oplus \mathbb{Z}$}
Given an abelian group, it is generally difficult to classify all
rings whose additive group is the group. However, when the abelian
group is small, one can work out of the classification by hand.
 Here, we classify all commutative $\mathbb{Z}$-algebras (i.e. rings)
 with identity whose
additive group is the free abelian group of rank two. The result is
used in Example \ref{ring_exa}.  Let $\A$ be such a ring. Its
additive group $(\A,+)$ is generated by 1 and $x$. We have $1*1=1,
1*x=x*1=x$, and $x*x=a1+bx$ where $a, b$ are arbitrary integers.
Obviously the ring structure of $\A$ is completely determined by
$(a,b)$. Let $\A'$ be another ring whose additive group is generated
by $1'$ and $x'$ with $x'*x'=a'1'+b'x'$. We have

\begin{proposition} The two rings $\A$ and $\A'$ are isomorphic if and only if
$b^2+4a=b'^2+4a' $ and $b\equiv b' \pmod{2})$. In other words, the
isomorphism type of $\A$ is completely determined by $(b^2+4a,$ $ b
\pmod{2})\in (\Z,\Z_2)$.
\end{proposition}

\noindent {\em Proof.} Suppose $\A\cong \A'$ as modules, and let $f:
\A \rightarrow \A'$ be an isomorphism. Then $f(1)=1', f(x)=k1'+lx'$
where $k, l$ are integers. Since $\{ f(1), f(x)\}$ spans the $\A'$
as an abelian group, we have $l=\pm 1$.

Since $f$ preserves the multiplication,
$f(x^2)=(k1'+lx')^2=k^21'+2klx'+l^2x'^2=k^21'+2klx'+l^2(a'1'+b'x')
=(k^2+l^2a')1'+(2kl+l^2b')x'$. We also have
$f(a1+bx)=a1'+b(k1'+lx')=(a+kb)1'+lbx'$. Since $x^2=a1+bx, $ we have

\[ \left\{ \begin{array}{ll}  k^2+l^2a'=a+kb & \hspace{3mm} (1) \\
2kl+l^2b'=lb& \hspace{3mm} (2) \end{array} \right. \]

The above two equations are equivalent to the following, where
(3) is obtained by taking (1)$-\frac{k}{l}*$(2), and (4) is $\frac{1}{l}*$(2):

\[ \left\{ \begin{array}{ll}  l^2a'-klb'-k^2=a & \hspace{3mm} (3) \\
2k+lb'=b& \hspace{3mm} (4) \end{array} \right. \]

Take $4*$(3)$+(4)^2$, we obtain:

\[ l^2(b'^2+4a')=b^2+4a \]

This is equivalent to $b'^2+4a'=b^2+4a$ since $l=\pm 1$.

The relation $b\equiv b' \pmod{2}$ follows from equation (4) by taking mod 2 in both sides.

Conversely, if $(a,b)$ and $(a',b')$ satisfy the two relations, we
let $l=1$, and $k=\frac{1}{2}(b-b')$. A straight forward computation
shows that equations (1) and (2) hold. This implies that the map $f:
A\rightarrow \A'$ defined by $f(1)=1', f(x)=k1'+x'$ is a ring
isomorphism.

\noindent
\begin{scriptsize}Department of Mathematics, George Washington
University, Washington, DC 20052. Email: lhelmeg@gwu.edu. \\
Department of Mathematics, George Washington University, Washington,
DC 20052.  Email: rong@gwu.edu
\end{scriptsize}

\end{document}